# ON THE CONCENTRATION OF MEASURE PHENOMENON FOR STABLE AND RELATED RANDOM VECTORS

By Christian Houdré and Philippe Marchal

*Université Paris XII and Georgia Institute of Technology, and Ecole Normale Supérieure*

Concentration of measure is studied, and obtained, for stable and related random vectors.

Let $X$ be a standard normal vector in $\mathbb{R}^d$ and let $f:\mathbb{R}^d \to \mathbb{R}$ be Lipschitz (with constant one) with respect to the Euclidean distance. A seminal result of Borell [3] and of Sudakov and Tsirel'son [7] asserts that for all $x > 0$,

$$(1) \qquad P(f(X) - m(f(X)) \geq x) \leq 1 - \Phi(x) \leq \frac{e^{-x^2/2}}{2},$$

where $m(f(X))$ is a median of $f(X)$ and where $\Phi$ is the (one-dimensional) standard normal distribution function. The inequality (1) has seen many extensions and to date, most of the conditions under which these developments hold require the existence of finite exponential moments for the underlying vector $X$. It is thus natural to explore the robustness of this "concentration phenomenon" and to study the corresponding results for stable vectors. It is the purpose of these notes to initiate this study and to present a few concentration results for stable and related vectors, freeing us from the exponential moment requirement. Our main result will imply that if $X$ is an $\alpha$-stable random vector in $\mathbb{R}^d$, then for all $x > 0$,

$$(2) \qquad P(f(X) - m(f(X)) \geq x) \leq 1 \wedge \frac{C(\alpha, d)}{x^\alpha},$$

where the constant $C(\alpha, d)$ will be explicit.

Let $X \sim ID(b, 0, \nu)$, that is, let $X$ be a $d$-dimensional infinitely divisible vector without Gaussian component. For all $u \in \mathbb{R}^d$, its characteristic function $\varphi_X$ is given by $\varphi_X(u) = e^{\psi(u)}$, with

$$\psi(u) = i\langle u, b\rangle + \int_{\mathbb{R}^d} (e^{i\langle u, y\rangle} - 1 - i\langle u, y\rangle \mathbf{1}_{\|y\|\leq 1})\nu(dy), \tag{3}$$

where $b \in \mathbb{R}^d$ and where $\nu \not\equiv 0$ (the Lévy measure) is a positive Borel measure without atom at the origin and such that $\int_{\mathbb{R}^d} (\|y\|^2 \wedge 1)\nu(dy) < +\infty$ (throughout, $\langle \cdot, \cdot \rangle$ and $\|\cdot\|$ are, respectively, the Euclidean inner product and norm in $\mathbb{R}^d$).

It is well known that there is a one-to-one relationship between Lévy processes and infinitely divisible laws. More precisely, if $\{X(t): t \geq 0\}$ is a Lévy process (without Gaussian component) on $\mathbb{R}^d$, then for all $t \geq 0$, $u \in \mathbb{R}^d$,

$$\varphi_{X(t)}(u) = Ee^{i\langle u, X(t)\rangle} = e^{t\psi(u)}, \tag{4}$$

where $\psi$ is as in (3). Hence, an infinitely divisible vector $X$ can be viewed as $X(1)$, the Lévy process $\{X(t): t \geq 0\}$ at time 1.

Recall also that $X$ is $\alpha$-stable, $0 < \alpha < 2$, if the measure $\nu$ is given, for any Borel set $B \in \mathcal{B}(\mathbb{R}^d)$, by

$$\nu(B) = \int_{S^{d-1}} \lambda(d\xi) \int_0^{+\infty} \mathbf{1}_B(r\xi) \frac{dr}{r^{1+\alpha}}, \tag{5}$$

where $\lambda$ is a finite positive measure on $S^{d-1}$, the unit sphere of $\mathbb{R}^d$, called the spherical component of the Lévy measure. $X$ is symmetric $\alpha$-stable ($S\alpha S$) if and only if $\lambda$ is symmetric, in which case,

$$\varphi_X(u) = \exp\left\{-c_\alpha \int_{S^{d-1}} |\langle u, \xi\rangle|^\alpha \lambda(d\xi)\right\},$$

where $c_\alpha = \frac{\sqrt{\pi}\Gamma((2-\alpha)/2)}{\alpha 2^\alpha \Gamma((1+\alpha)/2)}$. Moreover, $X$ is rotationally invariant if and only if $\lambda$ is uniform on $S^{d-1}$ and then

$$\varphi_X(u) = e^{-c_{\alpha,d}\|u\|^\alpha},$$

where $c_{\alpha,d} = c_\alpha \int_{S^{d-1}} |\langle u/\|u\|, \xi\rangle|^\alpha \lambda(d\xi)$ does not depend on $u \in \mathbb{R}^d$. In particular, if $\lambda$ is the uniform probability measure on $S^{d-1}$, $c_{\alpha,d} = \frac{\Gamma(d/2)\Gamma((2-\alpha)/2)}{\alpha 2^\alpha \Gamma((d+\alpha)/2)}$.

(We refer the reader to Sato's book [6] for a good introduction to Lévy processes and infinitely divisible laws.)

In order to prove our first theorem, we need the lemma below. For the mean rather than a median (and $x$ rather than $x/2$), the result is obtained in [4]. However, it is standard that applying this result to the function $g(y) = \min(d(y, A), x)$, $y \in \mathbb{R}^d$, where $A = \{f \leq m\}$ and $m$ is a median of $f$, leads to deviation from a median. Indeed, $Eg \leq x/2$, and therefore $g - Eg \geq x/2$ whenever $f - m \geq x$.

LEMMA 1. *Let $X \sim ID(b, 0, \nu)$ with $\nu$ boundedly supported, let $V^2 = \int_{\mathbb{R}^d} \|x\|^2 \nu(dx)$ and let $R = \inf\{\rho > 0 : \nu(\{x : \|x\| > \rho\}) = 0\}$. Then for any Lipschitz function (with constant 1) $f : \mathbb{R}^d \to \mathbb{R}$,*

$$(6) \quad P(f(X) - m(f(X)) \geq x) \leq \exp\left\{\frac{x}{2R} - \left(\frac{x}{2R} + \frac{V^2}{R^2}\right) \log\left(1 + \frac{Rx}{2V^2}\right)\right\},$$

*for all $x > 0$, and where $m(f(X))$ is a median of $f(X)$.*

Above (and below), the Lipschitz property is usually taken with respect to the Euclidean norm, that is, $f$ is Lipschitz if $\|f\|_{\text{Lip}} = \sup_{x \neq y} \frac{|f(x) - f(y)|}{\|x-y\|} < +\infty$; however, other norms could be considered [e.g., see Remark 2(ii)].

We can now state

THEOREM 1. *Let $X$ be an $\alpha$-stable vector with Lévy measure $\nu$ given by (5). Let $f : \mathbb{R}^d \to \mathbb{R}$ be such that $\|f\|_{\text{Lip}} \leq 1$. Then*

$$(7) \quad P(f(X) - m(f(X)) \geq x) \leq \frac{K_1 \lambda(S^{d-1})}{\alpha(2-\alpha)x^\alpha},$$

*whenever*

$$x \geq \left(\frac{K_2 \lambda(S^{d-1})}{\alpha(2-\alpha)}\right)^{1/\alpha},$$

*and where $m(f(X))$ is a median of $f(X)$ while $K_1$, $K_2$ are two absolute constants.*

PROOF. For any $R > 0$, we have the identity in distribution $X \stackrel{d}{=} Y^{(R)} + Z^{(R)}$, where $Y^{(R)}$ and $Z^{(R)}$ are mutually independent infinitely divisible vectors with respective characteristic function $\varphi_{Y^{(R)}} = e^{\psi_Y^{(R)}}$ and $\varphi_{Z^{(R)}} = e^{\psi_Z^{(R)}}$. For $u \in \mathbb{R}^d$, the exponents are given by

$$\psi_Z^{(R)}(u) = \int_{\|y\|>R} (e^{i\langle u,y\rangle} - 1) \nu_X(dy),$$

$$\psi_Y^{(R)}(u) = i\langle u, \tilde{b}\rangle + \int_{\|y\|\leq R} (e^{i\langle u,y\rangle} - 1 - i\langle u,y\rangle \mathbf{1}_{\|y\|\leq 1}) \nu_X(dy),$$

with

$$\tilde{b} = b - \int_{\|y\|>R} y \mathbf{1}_{\|y\|\leq 1} \nu_X(dy),$$

where the last integral is understood coordinatewise (and so is the above difference) and where $\nu_X$ is the Lévy measure of $X$.

Next,

$$(8) \quad P(f(X) - m(f(X)) \geq x) \leq P(f(Y^{(R)}) - m(f(X)) \geq x) + P(Z^{(R)} \neq 0).$$

Let us first estimate the second probability in (8) involving the compound Poisson random vector $Z^{(R)}$:

$$P(Z^{(R)} \neq 0) = 1 - P(Z^{(R)} = 0)$$

$$\leq 1 - \exp\left(-\int_{\|x\|>R} \nu_X(dx)\right)$$

$$= 1 - \exp\left(-\int_{S^{d-1}} \lambda(d\xi) \int_{\|r\xi\|>R} \frac{dr}{r^{1+\alpha}}\right)$$

(9)
$$= 1 - \exp\left(-\frac{\lambda(S^{d-1})}{\alpha} R^{-\alpha}\right)$$

$$= 1 - \exp(-C_2(\alpha,\lambda) R^{-\alpha})$$

$$\leq \frac{C_2}{R^\alpha},$$

where $C_2 := C_2(\alpha, \lambda) = \frac{\lambda(S^{d-1})}{\alpha}$.

Turning our attention to $Y^{(R)}$, we first compute the quantities involved in Lemma 1:

(10)
$$P(f(Y^{(R)}) - m(f(Y^{(R)})) \geq x)$$
$$\leq \exp\left\{\frac{x}{2R} - \left(\frac{x}{2R} + \frac{V^2}{R^2}\right)\log\left(1 + \frac{Rx}{2V^2}\right)\right\},$$

where

(11)
$$V^2 = \int_{\|x\|\leq R} \|x\|^2 \nu_X(dx)$$
$$= \int_{S^{d-1}} \lambda(d\xi) \int_{\|r\xi\|\leq R} \|r\xi\|^2 \frac{dr}{r^{1+\alpha}}$$
$$= \int_{S^{d-1}} \lambda(d\xi) \int_0^R r^2 \frac{dr}{r^{1+\alpha}}$$
$$= C_1(\alpha,\lambda) R^{2-\alpha},$$

with $C_1(\alpha,\lambda) = \frac{\lambda(S^{d-1})}{2-\alpha}$. Hence (10) becomes

$$P(f(Y^{(R)}) - m(f(Y^{(R)})) \geq x)$$

(12)
$$\leq \exp\left\{\frac{x}{2R} - \left(\frac{x}{2R} + \frac{C_1}{R^\alpha}\right)\log\left(1 + \frac{R^\alpha x}{2RC_1}\right)\right\}$$
$$:= H^{(R)}(x)$$

(13)
$$\leq \frac{e^{x/2R}}{(1 + R^\alpha x/(2RC_1))^{x/2R}},$$

where $C_1 := C_1(\alpha, \lambda) = \frac{\lambda(S^{d-1})}{2-\alpha}$.

Now rewrite the first probability in (8) as

$$
\begin{aligned}
(14) \quad & P(f(Y^{(R)}) - m(f(X)) \geq x) \\
& = P(f(Y^{(R)}) - m(f(Y^{(R)})) \geq x + m(f(X)) - m(f(Y^{(R)}))).
\end{aligned}
$$

We want to bound $|m(f(X)) - m(f(Y^{(R)}))|$ and, as it will become clear from the proof, only the case $m(f(X)) < m(f(Y^{(R)}))$ (which we assume) presents some interest and needs some work. To this end, remark that, for any $x \geq 0$ and any function $f$, we have

$$|P(f(X) \leq x) - P(f(Y^{(R)}) \leq x)| \leq P(X \neq Y^{(R)}) = P(Z^{(R)} \neq 0).$$

Set $P_m = P(f(X) \leq m(f(X))) \geq 1/2$. Then,

$$
\begin{aligned}
P_m - P(f(Y^{(R)}) & \leq m(f(X))) \\
& = P(f(X) \leq m(f(X))) - P(f(Y^{(R)}) \leq m(f(X))) \\
& \leq P(Z^{(R)} \neq 0).
\end{aligned}
$$

Moreover, if $f$ is Lipschitz with $\|f\|_{\text{Lip}} \leq 1$,

$$
\begin{aligned}
P_m - P(Z^{(R)} \neq 0) & \leq P(f(Y^{(R)}) \leq m(f(X))) \\
& = P(f(Y^{(R)}) - m(f(Y^{(R)})) \leq m(f(X)) - m(f(Y^{(R)}))) \\
& \leq H^{(R)}(m(f(Y^{(R)})) - m(f(X))),
\end{aligned}
$$

where the second inequality follows from Lemma 1 applied to $-f$, and with $H^{(R)}$ given in (12). Since $H^{(R)}$ is decreasing, set $I^{(R)}(y) = \sup\{z \geq 0, H^{(R)}(z) \geq y\}$. Thus, provided $P_m > P(Z^{(R)} \neq 0)$,

$$m(f(Y^{(R)})) - m(f(X)) \leq I^{(R)}(P_m - P(Z^{(R)} \neq 0)).$$

Choose $\delta \in (0, 1/2)$. Then, for every $R$ such that

$$(15) \quad R \geq (C_2/\delta)^{1/\alpha},$$

we have, using the same estimates as in (9), $P_m - P(Z^{(R)} \neq 0) \geq 1/2 - \delta$. Moreover, for every positive $A$, (13) entails

$$P(f(Y^{(R)}) - m(f(Y^{(R)})) \geq AR) \leq H^{(R)}(AR) \leq e^{A/2}\left(\frac{2C_1}{AR^\alpha}\right)^{A/2}.$$

Thus if

$$(16) \quad R \geq \left(\left(\frac{2C_1}{A}\right)^{A/2} \frac{e^{A/2}}{1/2 - \delta}\right)^{2/\alpha A},$$

then
$$H^{(R)}(AR) \leq 1/2 - \delta,$$

or, equivalently, $I^{(R)}(1/2 - \delta) \leq AR$. As a consequence, if $R$ satisfies both conditions (15) and (16), we have

$$|m(f(Y^{(R)})) - m(f(X))| \leq I^{(R)}(P_m - P(Z^{(R)} \neq 0)) \leq I^{(R)}(1/2 - \delta) \leq AR.$$

Using this together with (12) and (14) yields

$$P(f(Y^{(R)}) - m(f(X)) \geq (2+A)R)$$
$$\leq P(f(Y^{(R)}) - m(f(Y^{(R)})) \geq 2R) \leq \frac{eC_1}{R^\alpha}.$$

Setting $x = (2+A)R$, we obtain

(17) $$P(f(Y^{(R)}) - m(f(X)) \geq x) \leq \frac{eC_1(2+A)^\alpha}{x^\alpha}.$$

Finally, combining (17) and (9), we conclude that, for any $A > 0$,

(18) $$P(f(X) - m(f(X)) \geq x) \leq \frac{(eC_1 + C_2)(2+A)^\alpha}{x^\alpha},$$

whenever $x$ is large enough, so that there exists $\delta \in (0, 1/2)$ satisfying

(19) $$\frac{x}{2+A} \geq \left(\frac{C_2}{\delta}\right)^{1/\alpha},$$

and

(20) $$\frac{x}{2+A} \geq \left(\frac{2C_1}{A}\right)^{1/\alpha} \left(\frac{e^{A/2}}{1/2 - \delta}\right)^{2/\alpha A}.$$

For a given $A$, the domain of validity of (18) can be found by optimizing $\delta$ in (19) and (20). Taking, for instance, $A = 2$ leads [by equating the right-hand sides of (19) and (20)] to $\delta = C_2/2(eC_1 + C_2) = (2 - \alpha)/2(2 - \alpha + e\alpha)$, and so

$$P(f(X) - m(f(X)) \geq x) \leq \frac{4^\alpha(eC_1 + C_2)}{x^\alpha} = \frac{4^\alpha(2 - \alpha + e\alpha)\lambda(S^{d-1})}{\alpha(2 - \alpha)x^\alpha},$$

whenever

$$x \geq 4\left(\frac{2(2 - \alpha + e\alpha)\lambda(S^{d-1})}{\alpha(2 - \alpha)}\right)^{1/\alpha}. \qquad \square$$

REMARK 1. (i) The estimate in (7) is sharp in $x$, as can be seen by taking $X \sim S\alpha S$ a one-dimensional symmetric $\alpha$-stable random variable with parameter $\sigma > 0$ and characteristic function $\varphi_X(u) = e^{-\sigma^\alpha |u|^\alpha}$. In that case (e.g., see Proposition 1.2.15 in [5])

$$\lim_{x \to +\infty} x^\alpha P(X \geq x) = \sigma^\alpha A_\alpha,$$

where

$$A_\alpha = \begin{cases} \dfrac{1-\alpha}{2\Gamma(2-\alpha)\cos(\pi\alpha/2)}, & \alpha \neq 1, \\ \dfrac{1}{\pi}, & \alpha = 1, \end{cases}$$

and $\sigma^\alpha = 2\lambda(1)c_\alpha = \frac{\sqrt{\pi}\Gamma((2-\alpha)/2)}{2^\alpha \alpha \Gamma((1+\alpha)/2)} 2\lambda(1)$. For $d=1$, and $X$ symmetric, our constants are $C_1 = \frac{\lambda(1)+\lambda(-1)}{2-\alpha} = \frac{2\lambda(1)}{2-\alpha}$ and $C_2 = \frac{2\lambda(1)}{\alpha}$. Thus, the dependency in $\alpha$ in the constants of (7) is sharp as $\alpha \to 0$, but explodes as $\alpha \to 2$ (in contrast to $\sigma^\alpha A_\alpha$). This problem will be addressed in the sequel. We also note that the dependency on the dimension $d$ is sharp. Indeed, if $X$ is a stable vector in $\mathbb{R}^d$, then by a result of Araujo and Giné [1] (see, e.g., Theorem 4.4.8 in [5])

(21) $$\lim_{x \to +\infty} x^\alpha P(\|X\| \geq x) = c_\alpha \lambda(S^{d-1}) A_\alpha.$$

(ii) As usual left tails inequalities also follow from (7) by applying the result to $-f$ and, as is also classical, the estimates can equivalently be given in terms of enlargements of sets. For $\alpha > 1$, a median can be replaced by the mean (changing the range of $x$, too) by properly modifying the above proof or by using

$$E|f(X) - m(f(X))|$$
$$\leq 2\left(\frac{K_2 \lambda(S^{d-1})}{\alpha(2-\alpha)}\right)^{1/\alpha} + \frac{2}{\alpha-1} \frac{K_1 \lambda(S^{d-1})}{\alpha(2-\alpha)} \left(\frac{K_2 \lambda(S^{d-1})}{\alpha(2-\alpha)}\right)^{(1-\alpha)/\alpha},$$

which follows from integrating the tail inequality (7).

(iii) The methodology of proof presented above works as well for any infinitely divisible vector $X$. However, it requires estimates on $\int_{\|x\|>R} \nu_X(dx)$ and on $\int_{\|x\|\leq R} \|x\|^2 \nu_X(dx)$ which, in general, are unavailable in the absence of further knowledge of the Lévy measure. If the Lévy measure of $X$ has the form

(22) $$\nu'(B) = \int_{S^{d-1}} \lambda(d\xi) \int_0^{+\infty} \mathbf{1}_B(r\xi) \frac{L(r)\, dr}{r^{1+\alpha}},$$

for some slowly varying function $L$ on $[0, \infty)$, in which case $X$ is in the domain of attraction of a stable random vector with Lévy measure given by (5), the proof of Theorem 1 and standard estimates on regularly varying functions (see, e.g., [2], Theorem 1.5.11) give the following bound:

$$P(f(X) - m(f(X)) \geq x) \leq \frac{K_1 \lambda(S^{d-1}) L(x)}{\alpha(2-\alpha) x^\alpha}, \tag{23}$$

for every $x$ such that $L$ is locally bounded on $[x, \infty)$ and such that

$$\frac{x^\alpha}{L(x)} \geq K_2^\alpha,$$

where the constants $K_1$, $K_2$ are the same as in Theorem 1. A similar result can also be obtained when $X$ is in the domain of attraction of a stable vector with Lévy measure $\nu$. In that case, $L(r)$ in (22) should be replaced by $L(r, \xi)$; thus if $L_1(r) \leq L(r, \xi) \leq L_2(r)$, in (23), $L(x)$ should be replaced by $L_2(x)$.

When $\alpha$ is close to 2, the upper bound in Theorem 1 has the form $K\lambda(S^{d-1})/(2-\alpha)x^\alpha$ as soon as $x^\alpha > K'\lambda(S^{d-1})/(2-\alpha)$. We would like to obtain a better bound, namely of the form $K''\lambda(S^{d-1})/x^\alpha$, at the price of potentially strengthening the condition on $x$. To do so, we begin by a result improving Lemma 1. The setting and the notation below are as in Lemma 1; in addition, let $W^3 = \int_{\mathbb{R}^d} \|x\|^3 \nu(dx)$.

LEMMA 2. *If $V^2/W^3 > 2/R$, let $s_0$ be the unique positive solution of*

$$\frac{e^{sR} - 1}{sR} = \frac{RV^2}{W^3} - 1, \qquad s > 0,$$

*and let*

$$x_0 = 2\left(V^2 - \frac{W^3}{R}\right) s_0.$$

*Let $f : \mathbb{R}^d \to \mathbb{R}$ be such that $\|f\|_{\text{Lip}} \leq 1$. Then, if $x \leq x_0$,*

$$P(f(X) - Ef(X) \geq x) \leq \exp\left(\frac{-x^2}{4(V^2 - W^3/R)}\right),$$

*while for $x \geq x_0$,*

$$P(f(X) - Ef(X) \geq x) \leq K \exp\left(\frac{x}{R} - \left(\frac{x}{R} + \frac{2W^3}{R^3}\right) \log\left(1 + \frac{R^2 x}{2W^3}\right)\right),$$

*with*

$$K = \exp\left(\frac{-x_0^2}{4(V^2 - W^3/R)}\right)$$
$$\times \exp\left(-\frac{x_0}{R} + \left(\frac{x_0}{R} + \frac{2W^3}{R^3}\right) \log\left(1 + \frac{R^2 x_0}{2W^3}\right)\right). \tag{24}$$

PROOF. Recall that Theorem 1 in [4] asserts that

$$(25) \qquad P(f(X) - Ef(X) \geq x) \leq \exp\left\{-\int_0^x h^{-1}(s)\,ds\right\},$$

for all $0 < x < h(N^-)$, where $N = \sup\{t \geq 0 : Ee^{t\|X\|} < +\infty\}$ and where $h^{-1}$ is the inverse of $h(s) = \int_{\mathbb{R}^d} \|u\|(e^{s\|u\|} - 1)\nu(du)$, $0 < s < N$.

When the Lévy measure has its support in the Euclidean ball of center 0 and radius $R$ (in which case, $N = +\infty$), Lemma 1 follows by bounding the function $g_s(x) = e^{sx} - 1$ between 0 and $R$ by a linear interpolation, using also the convexity of the exponential. It is easily seen that, for every $x \in [0, R]$, the following improved inequality holds:

$$g_s(x) \leq sx + \frac{e^{sR} - 1 - sR}{R^2}x^2,$$

for all $s \geq 0$. Therefore,

$$h(s) \leq \left(V^2 - \frac{W^3}{R}\right)s + \frac{W^3}{R^2}(e^{sR} - 1)$$

$$\leq 2\max\left(\left(V^2 - \frac{W^3}{R}\right)s, \frac{W^3}{R^2}(e^{sR} - 1)\right).$$

So if $V^2/W^3 > 2/R$, there exists a unique positive $s_0$ such that

$$\frac{e^{s_0 R} - 1}{s_0 R} = \frac{RV^2}{W^3} - 1,$$

and for $s \leq s_0$,

$$h(s) \leq 2\left(V^2 - \frac{W^3}{R}\right)s,$$

while for $s \geq s_0$,

$$h(s) \leq 2\frac{W^3}{R^2}(e^{sR} - 1).$$

Let $x_0 = 2(V^2 - \frac{W^3}{R})s_0$. We have, for $t \leq x_0$,

$$h^{-1}(t) \geq \frac{t}{2(V^2 - W^3/R)},$$

while for $t \geq x_0$,

$$h^{-1}(t) \geq \frac{1}{R}\log\left(1 + \frac{R^2 t}{2W^3}\right).$$

The lemma follows. □

We are now ready to state our second theorem. We will express the deviation from the expectation here (since it exists). As already mentioned, a result in terms of the median can be easily derived.

THEOREM 2. *Using the notation of Theorem 1, assume $\alpha > 3/2$ and let $M = 1/(2-\alpha)$. Then there exists an absolute constant $K_3$ such that*

$$P(f(X) - Ef(X) \geq x) \leq \frac{K_3 \lambda(S^{d-1})}{x^\alpha}, \tag{26}$$

*for every $x$ satisfying*

$$x^\alpha \geq 4\lambda(S^{d-1}) M \log M \log(1 + 2M \log M).$$

PROOF. We use the notation of the proof of Theorem 1. The second and third moments of the Lévy measure $\nu_Y^{(R)}$ of $Y^{(R)}$ are given by

$$V^2 = \frac{\lambda(S^{d-1})}{2-\alpha} R^{2-\alpha},$$

and

$$W^3 = \frac{\lambda(S^{d-1})}{3-\alpha} R^{3-\alpha},$$

which entails

$$\frac{RV^2}{W^3} - 1 = M.$$

As $\alpha \geq 3/2$, $M \geq 2$ and

$$\frac{\log M}{R} \leq s_0 \leq 2\frac{\log M}{R},$$

and

$$\frac{2\lambda(S^{d-1}) M \log M}{(3-\alpha) R^{\alpha-1}} \leq x_0 \leq \frac{4\lambda(S^{d-1}) M \log M}{(3-\alpha) R^{\alpha-1}}. \tag{27}$$

So for $x \geq x_0$, Lemma 2 gives

$$P(f(Y^{(R)}) - Ef(Y^{(R)}) \geq x)$$
$$\leq K \exp\left(\frac{x}{R} - \left(\frac{x}{R} + \frac{2W^3}{R^3}\right) \log\left(1 + \frac{R^2 x}{2W^3}\right)\right).$$

From (24), (27) and since

$$\frac{2W^3}{R^3} \log\left(1 + \frac{R^2 x_0}{2W^3}\right) \leq \frac{x_0}{R},$$

it follows that

$$K \leq \exp\left(\frac{4\lambda(S^{d-1}) M \log M}{(3-\alpha) R^\alpha} \log(1 + 2M \log M)\right). \tag{28}$$

Suppose next that

(29) $$x^\alpha \geq 4\lambda(S^{d-1})M \log M \log(1 + 2M \log M).$$

Then setting $R = x$, it first follows (since $\alpha < 2$) that $K < e$, and second by (27) and since $M > 2$, that $x > x_0$. Therefore, since $W^3 = \lambda(S^{d-1})R^{3-\alpha}/(3-\alpha)$,

$$P(f(Y^{(R)}) - Ef(Y^{(R)}) \geq x)$$
$$\leq e \exp\left(1 - \left(1 + \frac{2\lambda(S^{d-1})}{(3-\alpha)x^\alpha}\right) \log\left(1 + \frac{(3-\alpha)x^\alpha}{2\lambda(S^{d-1})}\right)\right)$$
(30)
$$\leq e^2 \exp\left(-\log\left(1 + \frac{(3-\alpha)x^\alpha}{2\lambda(S^{d-1})}\right)\right)$$
$$\leq \frac{2e^2\lambda(S^{d-1})}{x^\alpha},$$

for all $x$ satisfying (29). Let us now estimate the difference between $Ef(X)$ and $Ef(Y^{(R)})$:

$$|Ef(X) - Ef(Y^{(R)})| = |E(f(Y^{(R)} + Z^{(R)}) - f(Y^{(R)}))\mathbf{1}_{\{Z^{(R)}\neq 0\}}|$$
$$\leq E\|Z^{(R)}\|$$
(31)
$$\leq \int_{\|x\|>R} \|x\|\nu_X(dx)$$
$$= \frac{\lambda(S^{d-1})}{\alpha - 1}R^{1-\alpha}$$
$$\leq \frac{x}{4\log 2\log(1 + 4\log 2)},$$

where we used the compound Poisson structure of $Z^{(R)}$ to get the next to last inequality, and our choice of $R = x$, $M > 2$, $\alpha - 1 > 1/2$, and the range of $x$ given by (29), to get the last one.

The estimate (31) as well as (30) and (9) finally give the result by proceeding as in the proof of Theorem 1 ($K_3 = 1 + 8e^2$ will do). □

REMARK 2. (i) Unless $\lambda(S^{d-1})$ is bounded above independently of $d$, Theorems 1 and 2 are not dimension-free, even for $X$ with independent components. This is to be expected in view of (21), and is in sharp contrast to the Gaussian case.

(ii) $X$ has independent components if and only if $\lambda$ is discrete and concentrated on the intersections of the axes of $\mathbb{R}^d$ with $S^{d-1}$. In that case, the natural Lipschitz property is with respect to the $\ell^1$-norm. Thus taking, in the Lévy–Khintchine representation (3), $S^{d-1}_{\|\cdot\|_1}$ (the $\ell^1$-unit ball of $\mathbb{R}^d$) instead of $S^{d-1}$, and correspondingly changing $\nu$ to $\nu_{\|\cdot\|_1}$, (7) continues to hold replacing $\lambda(S^{d-1})$ by $\lambda_{\|\cdot\|_1}(S^{d-1}_{\|\cdot\|_1})$, where $\lambda_{\|\cdot\|_1}$ is the "spherical component" of $\nu_{\|\cdot\|_1}$. In fact, for any norm of $\mathbb{R}^d$, a result similar to (7) continues to hold with the type of changes just described.

As already mentioned, the above theorems are not dimension-free, in particular, when the components of $X$ are independent and (for simplicity of notation) identically distributed. However, using Corollary 3 in [4], improved versions of Theorems 1 and 2 can be obtained with a mixture of "Lipschitz norms." More precisely, while we are not able to improve the constant in the upper bound of (26), the additional conditions we assume on the function $f$ enable us to extend (when $c$ below is small) the range on which our inequality holds. Again, for $X$ with i.i.d. components, the measure $\nu$ is concentrated on the axes of $\mathbb{R}^d$ (see [6], page 67), that is,

$$\nu(dx_1,\ldots,dx_d) = \sum_{k=1}^d \delta_0(dx_1)\cdots\delta_0(dx_{k-1})\tilde{\nu}(dx_k)\delta_0(dx_{k+1})\cdots\delta_0(dx_d).$$

Denoting by $e_1,\ldots,e_d$ the canonical basis of $\mathbb{R}^d$, Theorem 2 now becomes

THEOREM 3. *Let $X$ be a stable vector with index $\alpha > 3/2$ and Lévy measure $\nu$ given by* (31). *Let $f:\mathbb{R}^d \to \mathbb{R}$ be such that:*

(32) $$\sup_{x\in\mathbb{R}^d} \sum_{k=1}^d \int_\mathbb{R} |f(x+ue_k) - f(x)|^2 \tilde{\nu}(du) \leq a^2,$$

*and*

$$\sup_{x\in\mathbb{R}^d} |f(x+ue_k) - f(x)| \leq c|u|,$$

*for all $k=1,\ldots,d$, $u \in \mathbb{R}$. Then,*

$$P(f(X) - E(f(X)) \geq x) \leq \frac{K_4 \lambda(S^{d-1})c^\alpha + K_5 a^\alpha}{x^\alpha},$$

*whenever*

$$x^\alpha \geq 4(4\alpha)^{\alpha-1} c^{\alpha-1} \lambda(S^{d-1}),$$

*where $K_4$ and $K_5$ are absolute constants.*

Proof. We only briefly describe the changes to the previous proofs to get the result. First, and as previously done, decompose $X$ as $X \stackrel{d}{=} Y^{(R)} + Z^{(R)}$. Next, use Corollary 3 in [4]: under the assumptions on $f$ stated above

$$P(f(Y^{(R)}) - E(f(Y^{(R)})) \geq x) \leq \exp\left(-\frac{x}{4cR}\log\left(1 + \frac{cRx}{2a^2}\right)\right).$$

Then, proceed as in the proof of Theorem 1, using also the comparison between $Ef(X)$ and $Ef(Y^{(R)})$ given in the proof of Theorem 2: $|Ef(X) - Ef(Y^{(R)})| \leq \lambda(S^{d-1})R^{1-\alpha}/(\alpha-1)$. Hence, taking $R = Kx/(2\alpha c)$, for some constant $0 < K < 1$, chosen below, leads to:

$$P(f(X) - E(f(X)) \geq x) \leq \frac{(2\alpha)^\alpha \lambda(S^{d-1})c^\alpha + \alpha(4\alpha)^{\alpha/2}a^\alpha}{\alpha K^\alpha x^\alpha},$$

whenever

$$x^\alpha \geq \frac{(2\alpha)^{\alpha-1}}{\alpha-1}\frac{K^{1-\alpha}}{1-K}c^{\alpha-1}\lambda(S^{d-1}),$$

with $\alpha > 1$. The choices $K = 1/2$, $\alpha > 3/2$ yield the result. [A slightly improved result holds if in the integral in (32), $\mathbb{R}$ is replaced by $|u| \leq R$ where $R$ is chosen as above.] □

**Acknowledgments.** C. Houdré thanks the Laboratoire de Probabilités et Modèles Aléatoires, Université de Paris VI for its hospitality while part of this research was done. P. Marchal would like to thank the School of Mathematics at the Georgia Institute of Technology where part of this research was done.

## REFERENCES

[1] Araujo, A. and Giné, E. (1979). On tails and domains of attraction of stable measures in Banach spaces. *Trans. Amer. Math. Soc.* **248** 105–119. MR521695
[2] Bingham, N. H., Goldie, C. M. and Teugels, J. L. (1987). *Regular Variation.* Cambridge Univ. Press. MR898871
[3] Borell, C. (1975). The Brunn–Minkowski inequality in Gauss space. *Invent. Math.* **30** 207–216. MR399402
[4] Houdré, C. (2002). Remarks on deviation inequalities for functions of infinitely divisible random vectors. *Ann. Probab.* **30** 1223–1237. MR1920106
[5] Samorodnitsky, G. and Taqqu, M. (1994). *Stable Non-Gaussian Random Processes. Stochastic Models with Infinite Variance. Stochastic Modeling.* Chapman and Hall, New York. MR1280932
[6] Sato, K.-I. (1999). *Lévy Processes and Infinitely Divisible Distributions.* Cambridge Univ. Press. MR1739520
[7] Sudakov, V. N. and Tsirel'son, B. S. (1974). Extremal properties of half-spaces for spherically invariant measures. *Zap. Nauch. Sem. LOMI* **41** 14–24. [English translation *J. Soviet Math.* **9** (1978) 9–18.] MR365680

Laboratoire d'Analyse et
de Mathématiques Appliquées
CNRS UMR 8050
Universite Paris XII
Créteil 94010
France
and
School of Mathematics
Georgia Institute of Technology
Atlanta, Georgia
USA
e-mail: houdre@math.gatech.edu

DMA
Ecole Normale Supérieure
75005 Paris
France
e-mail: Philippe.Marchal@ens.fr